\documentclass{amsart}

\usepackage{hyperref}
\usepackage{hyperref}
\usepackage{amsmath,amsthm}
\usepackage{amssymb}
\usepackage{enumerate}
\usepackage{bbm}
\usepackage{cite}
\usepackage{hyperref}
\usepackage{graphicx}
\usepackage{mathtools}
\usepackage{enumerate}

\usepackage{bbm}
\usepackage[english]{babel}
\usepackage[utf8]{inputenc}
\usepackage{amsmath}
\usepackage{graphicx}
\usepackage[colorinlistoftodos]{todonotes}
\usepackage{amsmath, amsthm, amssymb}
\usepackage[bottom]{footmisc}
\usepackage{cite}
\usepackage{hyperref}
\usepackage{graphicx}

\usepackage{longtable}
\graphicspath{{figures/}} 
\usepackage{mathtools}

\newcommand{\Mod}[1]{\ (\mathrm{mod}\ #1)}

\DeclareMathOperator{\li}{li}

\usepackage{enumitem}

\usepackage{enumitem}

\usepackage{graphicx}

\usepackage[T1]{fontenc}

\newtheorem{theorem}{Theorem}[section]
\newtheorem{corollary}[theorem]{Corollary}

\theoremstyle{definition}

\newtheorem{remark}[theorem]{Remark}

\numberwithin{equation}{section}

\usepackage{amsmath}
\usepackage{multirow}
\usepackage{amsmath,amssymb, url,color}
\usepackage{mathtools}

\newcommand\halfopen[2]{\ensuremath{(#1,#2]}}

\makeatletter
\@namedef{subjclassname@2020}{\textup{2020} Mathematics Subject Classification}
\makeatother

\begin{document}


\title[On the second Hardy--Littlewood conjecture]{On the second Hardy--Littlewood conjecture}

\author[B. Chahal]{Bittu Chahal}
\address{Department of Mathematics, IIIT Delhi, New Delhi 110020, India}
\email{bittui@iiitd.ac.in}

\author[E. Elma]{Ertan Elma}
\address{Baku State University, Main Building, 3rd Floor, Mathematics Research Laboratory, Academician Zahid Khalilov Street 33, Baku, Azerbaijan, AZ 1148}
\email{ertan.elma@bsu.edu.az}

\author[N. Fellini]{Nic Fellini}
\address{Department of Mathematics and Statistics,
         Queen's University, Kingston, Canada,
         48 University Avenue, K7L 3N8}
\email{n.fellini@queensu.ca}

\author[A. Vatwani]{Akshaa Vatwani}
\address{Department of Mathematics, 
	Indian Institute of Technology Gandhinagar, Gandhinagar, 
	Gujarat 382355, India}
\email{akshaa.vatwani@iitgn.ac.in}

\author[D. N. T. Vo]{Do Nhat Tan Vo}
\address{University of Northern British Columbia, Department of Mathematics and Statistics, Prince George, BC V2N 4Z9 Canada}
\email{vo0@unbc.ca}

\maketitle


\begin{abstract}
The second Hardy--Littlewood conjecture asserts that the prime counting function $\pi(x)$ satisfies the subadditive inequality
	\begin{align*}
		\pi(x+y)\leqslant \pi(x)+\pi (y)
	\end{align*}
for all integers $x,y\geqslant 2$. 
By linking the subadditivity of $\pi(x)$ to the error term in the Prime Number Theorem, we obtain unconditional improvements on the range of $y$ for which $\pi(x)$ is known to be subadditive. Moreover, assuming the Riemann Hypothesis, we show that for all $\epsilon>0$, there exists $x_{\epsilon} \geqslant 2$ such that for all $x\geqslant x_\epsilon$ and $y$ in the range
\begin{align*}
		\frac{(2+\epsilon)\sqrt{x}\log^2x}{8\pi}\leqslant y\leqslant x,
\end{align*} 
the inequality $\pi(x+y)\leqslant \pi(x) + \pi(y)$ holds.
\end{abstract}

\section{Introduction}
Let $y\geqslant 2$ be an integer. Consider the number of primes in the first interval $\halfopen{0}{y}$ of length $y$ and compare it with the number of primes in another interval $\halfopen{x}{x+y}$ (of length $y$ again) for some integer $x\geqslant 2$. Which one has more primes?
In 1923, Hardy and Littlewood \cite{Hardy_Littlewood_1923} conjectured that the other intervals $\halfopen{x}{x+y}$ contain no more primes than the first one $\halfopen{0}{y}$, that is, 
\begin{align}\label{main inequality}
	\pi(x+y)-\pi(x)\leqslant  \pi (y)
\end{align}
for all integers $x,y\geqslant 2$. This is called the second Hardy--Littlewood conjecture. 

The first Hardy--Littlewood conjecture, on the other hand, is reserved for the well-known prime $k$-tuple conjecture whose special case is the twin prime conjecture. In 1973, Hensley and Richards \cite{Hensley_Richards_1973} established a striking connection between the $k$-tuple conjecture and the second Hardy--Littlewood conjecture. In particular, they showed that the two conjectures are incompatible, i.e., at most one of the conjectures is true. Moreover, assuming that the $k$-tuple conjecture is true, they showed that for sufficiently large values of $x$, there are infinitely many $y\geqslant 2$ for which the inequality (\ref{main inequality}) does not hold. Since the $k$-tuple conjecture is widely believed to be true, one expects that the inequality above is false for some integers $x,y\geqslant 2$ but no such counterexamples are known.


On the other hand, (\ref{main inequality}) is indeed true in several cases especially when the length $y$ of the interval is sufficiently large relative to the starting point $x$ of the other intervals being compared. Note that the subadditivity inequality $\pi(x+y)\leqslant \pi (x)+\pi(y)$ is symmetric in $x$ and $y$ and thus one can restrict attention to the case when $y\leqslant x$ without loss of generality. Landau \cite{Landau_1909} proved that 
\begin{align*}
	\pi(2x)\leqslant 2\pi (x)
\end{align*}   
for sufficiently large $x$ and Rosser and Schoenfeld \cite{Rosser_Schoenfeld_1975} showed that this is true for all $x\geqslant 2$. Karanikolov \cite{Karanikolov_1971} extended this result to $\pi(kx)< k\pi(x)$ for all $k\geqslant \sqrt{e}=1.648\dots$ and $x\geqslant 347 $ by using a result of Rosser, Schoenfeld and Yohe \cite{Rosser_Schoenfeld_Yohe_1968}. Panaitopol \cite{Panaitopol} proved that $\pi(kx)< k\pi(x)$ holds for fixed $k>1$ and sufficiently large $x$. By using the large sieve, Montgomery and Vaughan \cite{Montgomery_Vaughan_1973} proved that 
\begin{align*}
	\pi(x+y)\leqslant \pi (x)+2\pi(y)
\end{align*}
for all integers $x\geqslant 1,\, y\geqslant 2$. Udrescu \cite{Udrescu_1975} proved that if $\epsilon>0$, $x,y\geqslant 17$ and $x+y\geqslant1 +\exp\left(4+4/\epsilon\right)$, then
\begin{align*}
	\pi(x+y)\leqslant (1+\epsilon)\left(\pi(x)+\pi(y)\right)
\end{align*}
which indicates that the second Hardy--Littlewood conjecture could only be $\epsilon$-away from the truth. Moreover, Udrescu \cite{Udrescu_1975} showed that the original inequality
\begin{align}\label{subadditive_inequality}
	\pi(x+y)\leqslant \pi(x)+\pi(y)
\end{align}
holds as long as $\epsilon x\leqslant y\leqslant x$ for fixed $\epsilon\in\halfopen{0}{1}$ and sufficiently large $x$. Dusart \cite{Dusart_2002} improved on the range for $y$ for the validity of (\ref{subadditive_inequality}) to the range
\begin{align}\label{Dusart's_range}
	\frac{5x}{7\log x\log\log x} \leqslant y \leqslant x
\end{align}
 for $x\geqslant 5$. Here we would like to note that in \cite{Dusart_2002}, Dusart stated this range in a different form but as pointed out in Alkan's work \cite{Alkan_2022}, Dusart's range can be written as above by using the symmetry in the subadditivity inequality under consideration.    
 
 In this work we enlarge Dusart's range qualitatively depending on the shape of the error term in the Prime Number Theorem. Let 
 \begin{align*}
 	\li(x):=\int_{2}^{x}\frac{1}{\log u}\, du, \quad (x\geqslant 2)
 \end{align*}
be the logarithmic integral. The Prime Number Theorem is the assertion that 
\begin{align*}
	\pi(x)\sim \li (x)
\end{align*}
 as $x\rightarrow \infty$. Let $R(x)$ be a  function  such that there exist an absolute constant $C>0$ and a threshold $x_0\geqslant 2$ such that for all $x\geqslant x_0$, $R(x)$ is positive, nondecreasing  and  
\begin{align}\label{defn_R_x}
	\left|\pi(x)-\li(x)\right|\leqslant CR(x).
\end{align}
By a classical result of Littlewood \cite{Littlewood_1914}, \cite[Theorem 35, p. 103]{Ingham_1964}, we know that 
\begin{align*}
	\pi(x)-\li(x)=\Omega_{\pm}\left(\frac{x^{1/2}\log\log\log x}{\log x}\right)
\end{align*} 
as $x\rightarrow \infty$ where the notation $\Omega_{\pm}$ has the meaning that the limit superior (the limit inferior resp.) of the ratio of both sides is positive (negative resp.). Also the classical error term $O\left(x\exp\left(-c\sqrt{\log x}\right)\right)$ for some positive constant $c>0$ in the Prime Number Theorem saves an arbitrary power of $\log x$ and thus we may assume that the threshold $x_0$ is sufficiently large so that the function $R(x)$ in (\ref{defn_R_x}) satisfies
\begin{align*}
x^{1/2}\leqslant R(x)\leqslant \frac{x}{\log^3 x}
\end{align*}
for $x\geqslant x_0$.

\begin{theorem}\label{main_result}
	Let $x\geqslant   x_0$ where $x_0, C$ and $R(x)$ satisfy (\ref{defn_R_x}). If 
	\begin{align*}
			\frac{3CR(2x)\log^2 x}{\log\log x}\leqslant y\leqslant x,
	\end{align*}
	then the inequality
	\begin{align*}
		\pi(x+y)\leqslant \pi (x)+ \pi (y)
	\end{align*}
holds.
\end{theorem} 
 
The range for $y$ above improves on Dusart's range (\ref{Dusart's_range}) significantly. For some explicit choices for the function $R(x)$, Johnston and Yang \cite{Johnston_Yang_2023} proved that 
 \begin{align}\label{classical_error}
 	\left|\pi(x)-\li (x)\right|\leqslant 9.59x(\log x)^{0.515}\exp\left(-0.8274\sqrt{\log x}\right)
 \end{align}
for $x\geqslant 2$, and that 
\begin{align}\label{Vinogradov_Korobov_error}
	\left|\pi(x)-\li (x)\right|\leqslant 0.028x(\log x)^{0.801}\exp\left(-0.1853(\log x)^{3/5}(\log\log x)^{-1/5}\right)
\end{align}
for $x\geqslant 23$ and noted that the bound in (\ref{classical_error}) is superior to the one in (\ref{Vinogradov_Korobov_error}) for $\exp(59)\leqslant x \leqslant \exp\left(2.8\cdot10^{10}\right)$.  In \cite{Mossinghoff_Trudgian_Yang_2024}, Mossinghoff, Trudgian and Yang obtained that
\begin{align}\label{unconditional_bound}
		\left|\pi(x)-\li (x)\right|\ll  x\exp\left(-0.2123(\log x)^{3/5}(\log\log x)^{-1/5}\right)
\end{align}
without specifying the implied constant. These results are unconditional, but if the Riemann Hypothesis holds, then 
\begin{align}\label{conditional_bound}
	\left|\pi(x)-\li(x)\right|< \frac{1}{8\pi }\sqrt{x}\log x
\end{align}
for $x\geqslant 2657$ as shown by Schoenfeld \cite{Schoenfeld_1976}. By (\ref{conditional_bound}) and Theorem \ref{main_result}, the lower bound for $y$ becomes $\asymp x^{1/2}\log^3 x/ \log\log x$ under the Riemann Hypothesis. By the proof of Theorem \ref{main_result} below, we observe that if $R(x)\ll x^{\Theta}$ for some $\Theta<1$, then one can save a factor of $\log x /\log\log x$ in the lower bound for $y$. We reveal this fact under the Riemann Hypothesis in our next result.
\begin{corollary}\label{conditional_corollary_with_saving}
	Assume the Riemann Hypothesis. Then for all $\epsilon>0$, there exists $x_{\epsilon}\geqslant 2$ such that for all $x\geqslant x_\epsilon$ and 
	\begin{align*}
		\frac{(2+\epsilon)x^{1/2}\log^2 x}{8\pi}\leqslant y\leqslant x,
	\end{align*}
	the inequality $\pi(x+y)\leqslant \pi(x)+\pi (y)$ holds. 	
\end{corollary}

\begin{remark}\label{remark_r_1_r_2}
	The expression $2+\epsilon$ above can indeed be replaced by $(1+r_1(x))(2+r_2(x))$ where
	\begin{align*}
			r_1(x):&=\frac{2\log\left(0.08\log^{3}x\right)}{\log x^{1/2}}+\frac{35\log^2\left(0.08\log^{3}x\right)}{\log^2 x^{1/2}}=o_{x\rightarrow \infty}(1),
			\\r_2(x):&=-1+\left(1+\frac{0.08\log^3 x}{x^{1/2}}\right)^{1/2}\left(1+\frac{1}{\log x}\log\left(1+\frac{0.08\log^{3}x}{x^{1/2}}\right)\right) 
			\\&\nonumber\quad \quad  +\frac{\left(0.08\log^3 x\right)^{1/2}}{x^{1/4}}\left(\frac{1}{2}+\frac{\log\left(0.08\log^3 x\right)}{\log x}\right)=o_{x\rightarrow \infty}(1)
	\end{align*}
for $x\geqslant x_0':= 4\cdot 10^5$ which is a threshold up to where we verified the subadditivity inequality for all integers $x$ and $y$ with $2\leqslant y\leqslant x\leqslant x_0'$. The expression $(1+r_1(x))(2+r_2(x))$ tends to $2$ slowly as $(1+r_1(x_0'))(2+r_2(x_0'))=65.097\dots$.
\end{remark}

In the spirit of Corollary \ref{conditional_corollary_with_saving}, 
one may also ask whether the  partial verification of the Riemann Hypothesis up to some height could help to verify the subadditivity inequality (\ref{subadditive_inequality}).  More precisely, if $T_0$ is the largest known value such that the Riemann Hypothesis is true for all nontrivial zeros $\rho$ of $\zeta(s)$ with $\rm{Im}(\rho) \in (0,T_0]$, then it was recently shown by Johnston \cite{Johnston_2022} that the  bound  \eqref{conditional_bound} holds provided that
\begin{align} \label{johnston}
x \geqslant 2657 \: \text{ and }    \:
\frac{9.06}{\log\log x} \sqrt{\frac{x}{\log x}}  \leqslant T_0
\end{align}
which is an improvement on an earlier work of B\"uthe \cite{Buthe_2016}. 
By using this result, we obtain the following corollary. 
\begin{corollary} 
	\label{cor: RH up to height T}
Let $x\geqslant 4 \cdot 10^5$ and $r_1(x)$ and $r_2(x)$ be the functions in Remark \ref{remark_r_1_r_2}. Suppose that  the Riemann Hypothesis is true for all nontrivial zeros $\rho$ of $\zeta(s)$ with $\rm{Im}(\rho) \in (0,T_0]$. If
\begin{align*}
	\frac{(1+r_1(x))(2+r_2(x))x^{1/2}\log^2 x}{8\pi}\leqslant y\leqslant x
\end{align*}
and 
\begin{align*}
	\frac{9.06}{\log\log (x+y)} \sqrt{\frac{(x+y)}{\log (x+y)}}  \leqslant T_0,
\end{align*}
then the inequality $\pi(x+y)\leqslant \pi(x)+\pi(y)$ holds.
\end{corollary}
We remark here that the recent work by Platt and Trudgian \cite{Platt_Trudgian_2021} allows one to take $T_0 = 3 \cdot 10^{12}$.

As the inequality (\ref{subadditive_inequality}) holds for many pairs $x, y$, one may want to know an upper bound for the size of the exceptional set. Dusart \cite{Dusart_2002} and in a more general set-up Alkan \cite{Alkan_2022} considered this problem and proved that the size of the exceptional set is small by having a logarithmic saving. By Theorem \ref{main_result}, we improve on these results.
\begin{corollary}\label{corollary_exceptional_set}
	 We have 
	\begin{align*}
		\#\lbrace 2\leqslant y\leqslant x \leqslant X: \pi(x+y)> \pi (x)+\pi (y) \rbrace\ll  \frac{XR(2X)\log^2X}{\log\log X}.
	\end{align*}
	 In particular, the right-hand side above is
	 \begin{align*}
	 	\ll \frac{X^2\exp\left(-0.2123(\log 2X)^{3/5}(\log\log 2X)^{-1/5}\right)\log^2X}{\log\log X}
	 \end{align*}
unconditionally, and is
\begin{align*}
	\ll X^{3/2}\log^2X
\end{align*}
conditionally on the Riemann Hypothesis.
\end{corollary}

\begin{remark}
	Analogous results to Theorem \ref{main_result} and Corollaries \ref{conditional_corollary_with_saving} and \ref{corollary_exceptional_set} can be obtained in a more general setting. Let $P$ be a set of prime numbers and $\pi_P(x)$  be the number of primes $\leqslant x$ lying in the set $P$. In \cite[Theorem 1]{Alkan_2022}, Alkan proved that if 
	\begin{align*}
		\pi_{P}(x)=c\pi (x)+O\left(\frac{x}{\log^3 x}\right)
	\end{align*} 
for some constant $0<c\leqslant 1$, then for sufficiently large $x$ and $y$ satisfying 
\begin{align*}
	\frac{x}{\log x}\ll y\leqslant x,
\end{align*} 
the inequality $\pi_P(x+y)\leqslant \pi_P(x)+\pi_P(y)$ holds. A special case of this result gives the subadditivity inequality for the number of primes in a reduced residue class, \cite[Eq. (3.4)]{Alkan_2022}. Let $q\geqslant 2$ and $(a,q)=1$ and $\pi(x; q, a)$ be the number of prime numbers $p\leqslant x$ such that $p\equiv a \Mod{q}$. By the Prime Number Theorem for arithmetic progressions, \cite[Corollary 11.21]{Montgomery_Vaughan_2007_book}, we know that there exists a positive constant $c$ satisfying the following. For all $A>0$, there exists $x_0\geqslant 2$ and $C_A>0$ such that for all $x\geqslant x_0$, $q\leqslant (\log x)^A$ and $(a,q)=1$, we have 
\begin{align}\label{error_PNT_arithmetic_prog}
	\left|\pi(x; q, a)-\frac{\li(x)}{\varphi(q)}\right|\leqslant C_Ax\exp\left(-c\sqrt{\log x}\right)
\end{align} 
where $\varphi(.)$ is the Euler totient function. By following the proof of Theorem \ref{main_result}, one can unconditionally obtain that the inequality 
\begin{align}\label{ineq_arithmetic_prog}
	\pi(x+y; q, a)\leqslant \pi(x; q, a)+\pi(y; q, a), \quad \quad (q\leqslant (\log x)^A,\,  (a,q)=1)
\end{align}
holds for sufficiently large $x$ and $y$ satisfying 
\begin{align}\label{range_for_y_arithmetic_prog_unconditional}
\varphi(q)x\exp\left(-c\sqrt{\log x}\right)\frac{\log^2 x}{\log\log x}	\ll y\leqslant x
\end{align} 
where the implied constant depends only on $A$. Moreover, conditionally on the Generalized Riemann Hypothesis for all Dirichlet $L$-functions, the shape of the error term in (\ref{error_PNT_arithmetic_prog}) becomes $\ll x^{1/2}\log x$ for $q\leqslant x^{\frac{1}{2}-\epsilon}$ for arbitrarily small but fixed $\epsilon>0$ and thus the range for $y$ in (\ref{range_for_y_arithmetic_prog_unconditional}) for the validity of the inequality (\ref{ineq_arithmetic_prog}) becomes
\begin{align*}
	\varphi(q)x^{1/2}\log^2x \ll y \leqslant x
\end{align*}
by following the proof of Corollary \ref{conditional_corollary_with_saving}.
\end{remark}

\section{Proofs of Theorem \ref{main_result}  and Corollaries \ref{conditional_corollary_with_saving}, \ref{cor: RH up to height T} and \ref{corollary_exceptional_set}}
\subsection{Proof of Theorem \ref{main_result}}
Define 
\begin{align*}
	\Delta(x,y):=\pi(x)+\pi(y)-\pi(x+y)
\end{align*}
for $x,y\geqslant 2$. Thus, for the validity of (\ref{subadditive_inequality}), we need to show that $\Delta(x,y)\geqslant 0$. 

We can make two simplifications to the ranges of $x$ and $y$ considered. First, using Segal's criterion \cite{Segal_1962} for the validity of (\ref{subadditive_inequality}), we have verified (\ref{subadditive_inequality}) for all pairs $x$ and $y$ satisfying $x+y< 10^6$. Therefore, we may assume with no loss of generality that $x_0>4\cdot 10^5$. Secondly, by Dusart's result stated in (\ref{Dusart's_range}), we may assume, for $x\geqslant x_0$, that
\begin{align}\label{restriction1_on_y}
x^{1/2}\leqslant  y \leqslant \frac{x}{\log x }
\end{align}
since the remaining range $\frac{x}{\log x}<y\leqslant x$ is already covered in (\ref{Dusart's_range}). By the definition of the error term $R(x)$ in the approximation for $\pi(x)$ by $\li(x)$, we have 
\begin{align}\label{inequality_with_li}
	\nonumber \Delta(x,y)&=\li(x)+\li(y)-\li(x+y)
	 \\&\nonumber\quad +(\pi(x)-\li(x))+(\pi(y)-\li(y))-(\pi(x+y)-\li(x+y))
	\\&\geqslant \li(x)+\li(y)-\li(x+y)-3CR(2x)
\end{align} 
since $y\leqslant x$ and $R(x)$ is nondecreasing for $x\geqslant x_0$. We have 
\begin{align}\label{li_terms_first_lower_bound}
	 \nonumber \li(x)+\li(y)-\li(x+y)&=\int_{2}^{y}\frac{1}{\log u}\, du-\int_{x}^{x+y}\frac{1}{\log u}\, du
	 \\&\nonumber=\int_{2}^{y}\frac{1}{\log u}\, du-\int_{0}^{y}\frac{1}{\log (x+u)}\, du
	 \\&\nonumber=\int_{2}^{y}\frac{1}{\log u}\, du-\int_{2}^{y}\frac{1}{\log (x+u)}\, du-\int_{0}^{2}\frac{1}{\log (x+u)}\, du
	 \\& \geqslant \int_{2}^{y}\frac{1}{\log u}\left(1-\frac{\log u}{\log(x+u)}\right)\, du -\frac{2}{\log x}.
\end{align}
For $2\leqslant u \leqslant y \leqslant x/\log x$, we have 
\begin{align*}
1-\frac{\log u}{\log(x+u)} \geqslant 1 -\frac{\log(x/ \log x)}{\log (x+2)}  \geqslant  \frac{\log\log x}{\log x}. 
\end{align*}
Thus,
\begin{align*}
\li(x)+\li(y)-\li(x+y)\geqslant 	 \frac{\log\log x}{\log x} \int_{2}^{y}\frac{1}{\log u}\, du -\frac{2}{\log x}.
\end{align*}
On integration by parts, we have 
\begin{align}\label{li_integration_by_parts}
	\nonumber\int_{2}^{y}\frac{1}{\log u}\, du&=\frac{y}{\log y}-\frac{2}{\log 2}+\int_{2}^{y}\frac{1}{\log^2u}\, du
	\\&\geqslant 
	\frac{y}{\log x}-\frac{2}{\log 2}+\frac{x^{1/2}-2}{\log^2 x}
\end{align}
since $x^{1/2}\leqslant y\leqslant x$. Thus,
\begin{align}\label{li_terms_bound}
\nonumber	\li(x)+\li(y)-\li(x+y) &\geqslant \frac{\log\log x}{\log x}\left(\frac{y}{\log x}-\frac{2}{\log 2}+\frac{x^{1/2}-2}{\log^2 x} \right)-\frac{2}{\log x}
	\\&\geqslant \frac{y\log\log x}{\log^2 x}
\end{align}
since 
\begin{align*}
\frac{\log\log x}{\log x}\left(-\frac{2}{\log 2}+\frac{x^{1/2}-2}{\log^2x}\right)-\frac{2}{\log x}>0
\end{align*}
for $x\geqslant 4\cdot 10^{5}$. By (\ref{inequality_with_li}) and (\ref{li_terms_bound}), we have
\begin{align*}
	 \Delta(x,y)\geqslant \frac{y\log\log x}{\log^2 x}-3CR(2x)
\end{align*}
which is nonnegative if 
\begin{align*}
	y\geqslant \frac{3CR(2x)\log^2 x}{\log\log x}
\end{align*}
where $x\geqslant x_0$. This finishes the proof of Theorem \ref{main_result}. 

\subsection{Proof of Corollaries \ref{conditional_corollary_with_saving} and \ref{cor: RH up to height T}}
We first consider the conditional result, Corollary \ref{conditional_corollary_with_saving}. We apply the argument in the proof of Theorem \ref{main_result} above by modifying the range (\ref{restriction1_on_y}) on $y$ and using the Riemann Hypothesis.  
Let 
\begin{align*}
	c_1&=\frac{3\sqrt{ 2}}{8 \pi},
	\\x_0'&=4\cdot 10^5.
\end{align*}
On the Riemann Hypothesis, we already know that the subadditivity inequality (\ref{subadditive_inequality}) holds if 
\begin{align}\label{conditional_range_no_saving}
	\frac{c_1x^{1/2}\log(2x) \log^2 x}{\log\log x}\leqslant y \leqslant x, \quad (x\geqslant x_0')
\end{align}
by using  the conditional bound (\ref{conditional_bound}) and Theorem \ref{main_result}.  Thus, going over the proof of Theorem \ref{main_result} again, we may restrict $y$ to the range 
\begin{align}\label{conditional_restriction_on_y}
	x^{1/2}\leqslant y \leqslant 	\frac{c_1x^{1/2} \log(2x) \log^2 x}{\log\log x}\leqslant c_2x^{1/2}\log^3 x, \quad\quad  (x\geqslant x_0')
\end{align}
where 
\begin{align*}
	c_2=0.08.
\end{align*}
By (\ref{li_terms_first_lower_bound}), we have 
\begin{align*}
	 \li(x)+\li(y)-\li(x+y)\geqslant \int_{2}^{y}\frac{1}{\log u}\left(1-\frac{\log u}{\log(x+u)}\right)\, du -\frac{2}{\log x}.
\end{align*}
For $2\leqslant u \leqslant y \leqslant c_2x^{1/2}\log^3 x$, we have  
\begin{align*}
	1-\frac{\log u}{\log (x+u)}\geqslant1-\frac{\log\left(c_2x^{1/2}\log^3x\right)}{\log x}=\frac{1}{2}-\frac{\log\left(c_2\log^{3}x\right)}{\log x}>0
\end{align*}
for $x>x_0'$. Thus, similar to (\ref{li_integration_by_parts}) and (\ref{li_terms_bound}), we have 
\begin{align*}
 &\li(x)+\li(y)-\li(x+y)
 \\& \geqslant \left(\frac{1}{2}-\frac{\log\left(c_2\log^{3}x\right)}{\log x}\right)	\left(\frac{y}{\log \left(c_2x^{1/2}\log^3 x\right)}-\frac{2}{\log 2}+\frac{x^{1/2}-2}{\log^2 \left(c_2x^{1/2}\log^3 x\right)} \right)
 \\&\quad \quad \quad -\frac{2}{\log x}
 \\&\geqslant \left(\frac{1}{2}-\frac{\log\left(c_2\log^{3}x\right)}{\log x}\right)\frac{y}{\log \left(c_2x^{1/2}\log^3 x\right)}
\end{align*}
for $x\geqslant x_0'$ and $y$ in the range (\ref{conditional_restriction_on_y}). We have 
\begin{align*}
	\frac{1}{\log \left(c_2x^{1/2}\log^3 x\right)}=\frac{2}{\log x}\left(1-\frac{\log\left(c_2\log^3 x\right)}{\log \left(c_2x^{1/2}\log^3 x\right)}\right)\geqslant \frac{2}{\log x}\left(1-\frac{\log\left(c_2\log^{3}x\right)}{\log x^{1/2}}\right).
\end{align*}
Thus
\begin{align*}
\li(x)+\li(y)-\li(x+y)&\geqslant  \left(\frac{1}{2}-\frac{\log\left(c_2\log^{3}x\right)}{\log x}\right)\frac{2y}{\log x}\left(1-\frac{\log\left(c_2\log^{3}x\right)}{\log x^{1/2}}\right)
\\&=\frac{y}{\log x}\left(1-\frac{\log\left(c_2\log^{3}x\right)}{\log x^{1/2}}\right)^2.
\end{align*}
for $x\geqslant x_0'$ and $y$ in the range (\ref{conditional_restriction_on_y}).
Note that 
\begin{align*}
	\frac{\log\left(c_2\log^{3}x\right)}{\log x^{1/2}}<0.8
\end{align*}
for $x\geqslant x_0'$ and $(1-z)^{2}\geqslant (1+2z+35 z^2)^{-1}$ for $0<z<0.8$. Thus
\begin{align}\label{final_conditional_li_bound}
\li(x)+\li(y)-\li(x+y)&\geqslant 	\frac{y}{\log x}(1+r_1(x))^{-1}
\end{align}
where 
\begin{align*}
	r_1(x)=\frac{2\log\left(c_2\log^{3}x\right)}{\log x^{1/2}}+\frac{35\log^2\left(c_2\log^{3}x\right)}{\log^2 x^{1/2}}=o_{x\rightarrow \infty}(1).
\end{align*}

 Since $y$ satisfies (\ref{conditional_restriction_on_y}), one can also improve on the constant $3$ in (\ref{inequality_with_li}). Instead of $-3CR(2x)$ in (\ref{inequality_with_li}), the error we obtain is 
\begin{align}\label{conditional_error}
	\nonumber &\geqslant -\frac{1}{8\pi}x^{1/2}\log x -\frac{1}{8\pi}y^{1/2}\log y-\frac{1}{8\pi}\left(x+y\right)^{1/2}\log(x+y)
	\\&\geqslant -\frac{(2+r_2(x))}{8\pi}x^{1/2}\log x
\end{align}
for  $y$ in (\ref{conditional_restriction_on_y}),  where 
\begin{align*}
	r_2(x)&=-1+\left(1+\frac{c_2\log^3 x}{x^{1/2}}\right)^{1/2}\left(1+\frac{1}{\log x}\log\left(1+\frac{c_2\log^{3}x}{x^{1/2}}\right)\right) 
\\&\quad \quad  +\frac{\left(c_2\log^3 x\right)^{1/2}}{x^{1/4}}\left(\frac{1}{2}+\frac{\log\left(c_2\log^3 x\right)}{\log x}\right)=o_{x\rightarrow\infty}(1).      
\end{align*}
 By (\ref{final_conditional_li_bound}) and $\ref{conditional_error}$, we have 
\begin{align*}
\Delta(x,y)\geqslant 		\frac{y}{\log x}(1+r_1(x))^{-1}-\frac{(2+r_2(x))}{8\pi}x^{1/2}\log x
\end{align*}
which is nonnegative if 
\begin{align*}
	y\geqslant \frac{(1+r_1(x))(2+r_2(x))}{8\pi}x^{1/2}\log^2 x.
\end{align*}
This completes the proof of Corollary \ref{conditional_corollary_with_saving}.
The proof of Corollary \ref{cor: RH up to height T} is completely identical, since the  condition 
\begin{align*}
\frac{9.06}{\log\log (x+y)} \sqrt{\frac{(x+y)}{\log (x+y)}}  \leqslant T_0. 
\end{align*}
in the hypothesis of this corollary allows us to invoke Johnston's result (see \eqref{johnston}) and use the bound $|\pi(t)- \li(t)|< \frac{1}{8\pi} \sqrt t
 \log t$ for $t =x, y$ and $ x+y$. 
\subsection{Proof of Corollary \ref{corollary_exceptional_set}}
For Corollary \ref{corollary_exceptional_set}, we follow Alkan's  argument \cite{Alkan_2022} and observe that the number of exceptions to (\ref{subadditive_inequality}) with $2\leqslant y \leqslant x \leqslant X$ is at most
\begin{align}\label{two_double_sums}
	\sum_{x\leqslant X^{1/2}}\sum_{y\leqslant x}1+\sum_{X^{1/2}\leqslant x\leqslant X}\sum_{y\leqslant 	\frac{3CR(2x)\log^2x}{\log\log x}}1
\end{align}
by Theorem \ref{main_result}. The first term above is $\ll X $ and the second term
\begin{align*}
	\sum_{X^{1/2}\leqslant x\leqslant X}\sum_{y\leqslant 	\frac{3CR(2x)\log^2x}{\log\log x}}1&\ll\sum_{X^{1/2}\leqslant x\leqslant X} \frac{R(2x)\log^2x}{\log\log x}
	\\&\ll \frac{XR(2X)\log^2X}{\log\log X}
\end{align*}
which gives the first result in Corollary \ref{corollary_exceptional_set}. By using the unconditional bound (\ref{unconditional_bound}), the second assertion in Corollary \ref{corollary_exceptional_set} follows. For the third result in Corollary \ref{corollary_exceptional_set}, we replace the upper bound for $y$ in the second term in (\ref{two_double_sums}) by $\ll x^{1/2}\log^2 x$ and the desired result follows.  

\section*{Acknowledgments}
We are grateful to the Pacific Institute for the Mathematical Sciences (PIMS), Collaborative Research Groups (CRG) on $L$-functions,  Alia Hamieh, Ghaith Hiary, Habiba Kadiri,  Allysa Lumley, Greg Martin and Nathan Ng for organizing the summer school Inclusive Paths in Explicit Number Theory at the University of British Columbia, Okanagan, July 2023, where this project has started. We thank Tim Trudgian for suggesting this problem to study and providing us with valuable ideas and references. We are also grateful to  M. Ram Murty and Olivier Ramar\'e for their helpful comments and suggestions.   The second author was supported by a University of Lethbridge postdoctoral fellowship in 2023 and currently he is supported as a Junior Researcher at the Mathematics Research Lab at Baku State University.


\begin{thebibliography}{9}
	
\bibitem{Alkan_2022} E. Alkan, A generalization of the Hardy--Littlewood conjecture, Integers {\bf 22} (2022), Paper No. A53, 21 pp.; MR4441647.	
	
\bibitem{Buthe_2016} J. B\"uthe,  Estimating $\pi(x)$ and related functions under partial RH assumptions, Math. Comp. \textbf{85} (2016), no.~301, 2483--2498; MR3511289. 	
	
\bibitem{Dusart_2002} P. Dusart, Sur la conjecture $\pi(x+y)\le\pi(x)+\pi(y)$, Acta Arith. {\bf 102} (2002), no.~4, 295--308; MR1889398.
	
\bibitem{Hensley_Richards_1973}	D. Hensley and I. Richards, Primes in intervals, Acta Arith. {\bf 25} (1973/74), 375--391; MR0396440.
	
	
	
\bibitem{Hardy_Littlewood_1923} G.~H. Hardy and J.~E. Littlewood, Some problems of `Partitio numerorum'; III: On the expression of a number as a sum of primes, Acta Math. {\bf 44} (1923), no.~1, 1--70; MR1555183.

\bibitem{Ingham_1964} A.~E. Ingham, {\it The distribution of prime numbers}, Cambridge Tracts in Mathematics and Mathematical Physics, No. 30, Stechert-Hafner, Inc., New York, 1964; MR0184920.

\bibitem{Johnston_2022} D.~R. Johnston, Improving bounds on prime counting functions by partial verification of the Riemann hypothesis, Ramanujan J. {\bf 59} (2022), no.~4, 1307--1321; MR4507211. 

\bibitem{Johnston_Yang_2023} D.~R. Johnston and A. Yang, Some explicit estimates for the error term in the prime number theorem, J. Math. Anal. Appl. {\bf 527} (2023), no.~2, Paper No. 127460, 23 pp.; MR4601071.


\bibitem{Karanikolov_1971} H.~N. Karanikolov, On some properties of function $\pi (x)$, Univ. Beograd. Publ. Elektrotehn. Fak. Ser. Mat. Fiz. No. 357-380 (1971), 357, 29--30; MR0297719.

\bibitem{Landau_1909} E. Landau, {\it Handbuch der Lehre von der Verteilung der Primzahlen. 2 B\"ande}, Chelsea, New York, 1953; MR0068565.

\bibitem{Littlewood_1914} J. E. Littlewood, Sur la distribution des nombres premiers, C. R. Acad. Sci., Paris 158, (1914), 1869--1872.


\bibitem{Montgomery_Vaughan_1973} H.~L. Montgomery and R.~C. Vaughan, The large sieve, Mathematika {\bf 20} (1973), 119--134; MR0374060.

\bibitem{Montgomery_Vaughan_2007_book} H.~L. Montgomery and R.~C. Vaughan,\textit{ Multiplicative number theory I, Classical theory}, Cambridge Studies in Advanced Mathematics, Vol. 97, Cambridge University Press, Cambridge, 2007.

\bibitem{Mossinghoff_Trudgian_Yang_2024} M.~J. Mossinghoff, T.~S. Trudgian and A. Yang, Explicit zero-free regions for the Riemann zeta-function, Res. Number Theory {\bf 10} (2024), no.~1, Paper No. 11, 27 pp.; MR4688363.


\bibitem{Panaitopol} L. Panaitopol, Eine Eigenschaft der Funktion \"uber die Verteilung der Primzahlen, Bull. Math. Soc. Sci. Math. R. S. Roumanie (N.S.) {\bf 23(71)} (1979), no.~2, 189--194; MR0542926.

\bibitem{Platt_Trudgian_2021} D. Platt, T. Trudgian, The Riemann hypothesis is true up to $3\cdot 10^{12}$. Bull. Lond. Math. Soc. {\bf 53}(2021), no.~3, 792--797; MR4275089. 

\bibitem{Rosser_Schoenfeld_1975} J.~B. Rosser and L. Schoenfeld, Sharper bounds for the Chebyshev functions $\theta (x)$ and $\psi (x)$, Math. Comp. {\bf 29} (1975), 243--269; MR0457373.

\bibitem{Rosser_Schoenfeld_Yohe_1968} J.~B. Rosser, L. Schoenfeld and J.~M. Yohe, Rigorous computation and the zeros of the Riemann zeta-function. (With discussion), in {\it Information Processing 68 (Proc. IFIP Congress, Edinburgh, 1968), Vol. 1: Mathematics, Software}, pp. 70--76, North-Holland, Amsterdam; MR0258245.

\bibitem{Schoenfeld_1976} L. Schoenfeld, Sharper bounds for the Chebyshev functions $\theta (x)$ and $\psi (x)$. II, Math. Comp. {\bf 30} (1976), no.~134, 337--360; MR0457374.

\bibitem{Segal_1962} S. L. Segal, On $\pi (x+y)\leq \pi(x)+\pi (y)$, Trans. Amer. Math. Soc. {\bf 104} (1962), 523--527; MR0139586.

\bibitem{Udrescu_1975} V.~\c S. Udrescu, Some remarks concerning the conjecture $\pi (x+y)\leq \pi (x)+\pi (y)$, Rev. Roumaine Math. Pures Appl. {\bf 20} (1975), no.~10, 1201--1209; MR0392868.


\end{thebibliography}
\end{document}